\documentclass[smallextended]{svjour3}
\usepackage{latexsym}
\usepackage{amssymb}
\usepackage{epsfig}
\usepackage{amsmath}
\usepackage[mathscr]{eucal}
\usepackage{subfigure}
\usepackage{graphicx}

\newtheorem{Lemma}{Lemma}
\newtheorem{Proposition}{Proposition}
\newtheorem{Theorem}[Lemma]{Theorem}

\renewcommand{\qed}{\hfill{\ \ \rule{2mm}{2mm}} \vspace{0.2in}}

\renewcommand{\thefigure}{\arabic{figure}}
\begin{document}

\title{Law of Large Numbers for Permanents of Random Constrained Matrices}
\author{ \textbf{Ghurumuruhan Ganesan}
\thanks{E-Mail: \texttt{gganesan82@gmail.com} } \\
\ \\
Institute of Mathematical Sciences, HBNI, Chennai}
\date{}
\maketitle

\begin{abstract}
Permanents of random matrices with independent and identically distributed (i.i.d.) entries have extensively studied in literature and convergence and concentration properties are known under varying assumptions on the distributions. In this paper we consider constrained~\(n \times n\) random matrices where each row has a deterministic number~\(r = r(n)\) zero entries and the rest of the entries are independent random variables that are positive almost surely. The positions of the zeros within each row are also random and we establish sufficient conditions for the existence of a weak law of large numbers for the permanent, appropriately scaled and centred.  As a special case, we see that if~\(r\) grows faster than~\(\sqrt{n},\) then the permanent of a randomly chosen constrained~\(0-1\) matrix is concentrated around the van der Waerden bound with high probability, i.e. with probability converging to one as~\(n \rightarrow \infty.\)


\vspace{0.1in} \noindent \textbf{Key words:} Weak Law of Large Numbers, Row Constrained Random Matrices, van der Waerden bound, Second moment method.

\vspace{0.1in} \noindent \textbf{AMS 2000 Subject Classification:} Primary: 60C05, 05A16.
\end{abstract}

\bigskip

\renewcommand{\theequation}{\thesection.\arabic{equation}}
\setcounter{equation}{0}
\section{Introduction} \label{intro}


Estimating the permanents of random matrices is of great interest from both theoretical and application perspectives and  for matrices with independent and identically distributed (i.i.d.) entries, many convergence properties are known under varying assumptions on the distributions. For example, (Rempala and Gupta, 1999) have studied random matrices with independent entries and have obtained central limit theorems for the permanent as an extension of Girko's theorem for determinants. Later (Tao and Vu, 2009) used concentration techniques and showed that the permanent of a random~\(n \times n\) matrix with entries in~\(\{-1,1\}\) is non-zero with high probability, i.e. with probability converging to one as~\(n \rightarrow \infty\) and also estimated the typical value of the permanent. Recently, (Kwan and Sauermann, 2020) have obtained an analogue of the Tao-Vu result for the permanents of random symmetric matrices and for a recent survey of combinatorial properties of random matrices, we refer to (Vu, 2020).

In this paper we study permanents of~\(n \times n\) random matrices where each row has a deterministic number zero entries and the rest of the entries are i.i.d. random variables that are positive almost surely. The positions of the zeros within each row are also random and we use the second moment method to establish sufficient conditions on the mean and second moment that guarantee the existence of a weak law of large numbers for the permanent, appropriately scaled and centred.   In particular, for~\(0-1\) matrices, this also provides conditions under which the permanent of a randomly chosen constrained~\(0-1\) matrix is concentrated around the van der Waerden bound with high probability. (For completeness, we recall that the van der Waerden bound (Egorychev, 1981), (Falikman, 1981) asserts that the permanent of a deterministic~\(n \times n\) matrix containing  exactly~\(r\) ones in each row and column must be at least~\(r^{n}\frac{n!}{n^n}\)).




The paper is organized as follows. In the rest of the section, we state our main result Theorem~\ref{main_thm} regarding the~\(L^2-\)convergence of the permanent of constrained random matrices. Next in Section~\ref{prelim}, we collect preliminary estimates regarding the mean and the second moment of the permanent of~\(X\) and finally, in Section~\ref{main_sec}, we prove Theorem~\ref{main_thm}.

\subsection*{\em Constrained Random Matrices}
For integer~\(n\geq 1\) let~\(Z := [Z_{i,j}]_{1 \leq i,j \leq n}\) be a random matrix whose entries~\(Z_{i,j}, 1 \leq i,j \leq n\) are independent and identically distributed (i.i.d.) random variables that are positive a.s. with mean~\(\mathbb{E}Z_{1,1}\) and second moment~\(\delta = \delta(n) = EZ^2_{1,1}.\) The permanent of~\(Z\) is defined as
\begin{equation}\label{per_def}
per(Z) := \sum_{\sigma} \prod_{i=1}^{n}Z_{i,\sigma(i)},
\end{equation}
where the summation is over all permutations~\(\sigma\) of~\(\{1,2,\ldots,n\}.\)

We introduce constraints in the number of non-zero entries of~\(Z\) as follows. Let~\(1 \leq r_1,\ldots,r_n \leq n\) be deterministic integers and let~\({\cal M}\) be the set of all~\(n \times n\) matrices that contain exactly~\(r_i\) ones in the~\(i^{th}\) row and zeros in the rest of the positions, for each~\(1 \leq i \leq n.\) Let~\(X = [X_{i,j}]_{1 \leq i,j \leq n}\) be uniformly randomly chosen from~\({\cal M},\) independent of~\(Z\) and form the termwise product matrix~\(Y = [X_{i,j} \cdot Z_{i,j}]_{1 \leq i,j \leq n}.\) For~\(1 \leq i \leq n,\) the matrix~\(Y\) has exactly~\(r_i\) non-zero entries in row~\(i\) and we define the random matrix~\(Y\) on the probability space~\((\Omega,{\cal F}, \mathbb{P}).\)

The following result determines sufficient conditions for the~\(L^2-\)convergence of the permanent~\(T_n = per(Y)\) of the matrix~\(Y,\) appropriately scaled and centred. For two sequences~\(\{a_n\}_{n \geq 1}\) and~\(\{b_n\}_{n \geq 1},\) we say that~\(b_n = o(a_n)\) if~\(\frac{b_n}{a_n}\longrightarrow 0\) as~\(n \rightarrow \infty\) and~\(b_n = O(a_n)\) if there exists a constant~\(C > 0\) such that~\(b_n \leq C a_n.\) Throughout constants do not depend on~\(n.\)
\begin{Theorem}\label{main_thm}
Setting~\(\mu_n := \prod_{i=1}^{n} r_i \cdot \frac{\nu^{n} n!}{n^{n}},\) we have that~\(\mathbb{E}T_n = \mu_n.\) Moreover, if~\(r_{low} = \min_{1 \leq i \leq n} r_i \) and~\(r_{up} = \max_{1 \leq i \leq n} r_i\) are such that
\begin{equation}\label{r_cond}
\frac{\delta}{\nu^2 r_{low}} = o\left(\frac{1}{\sqrt{n}}\right)\text{ and }\frac{\delta }{\nu^2 r_{low}} - \frac{1}{r_{up}} = o\left(\frac{1}{n}\right)
\end{equation}
then
\begin{equation}\label{m_conv}
\mathbb{E}\left(\frac{T_n}{\mu_n}-1\right)^2 \longrightarrow 0 \text{ as } n\rightarrow \infty.
\end{equation}
\end{Theorem}
For illustration, we consider the example of~\(0-1\) matrices and set~\(Z_{i,j}  = 1\) for all~\(i,j.\) In this case~\(\delta = \nu^2 =1\) and the conditions in~(\ref{r_cond}) then hold if the matrices in~\({\cal M}\) are reasonably homogenous in terms of the number of ones in each row and is satisfied if for example
\[ (\sqrt{n}-1) \log{n} \leq r_i \leq (\sqrt{n}+1)\log{n},\] for each~\(1 \leq i \leq n.\) In the case of complete homogeneity, i.e. if~\(r_i = r = r(n)\) for each~\(1 \leq i \leq n,\) then~\(\mu_n = r^{n} \cdot \frac{n^{n}}{n!}\) is the van der Waerden bound and condition~(\ref{r_cond}) is equivalent to~\(\frac{r}{\sqrt{n}} \longrightarrow \infty\) as~\(n \rightarrow \infty.\)


To prove Theorem~\ref{main_thm}, we use the second moment method: We first obtain first and second moment bounds for the permanent of~\(Y\) and then show that if~(\ref{r_cond}) holds, then the permanent variance, appropriately scaled, converges to zero as~\(n \rightarrow \infty.\) This obtains the desired convergence~(\ref{m_conv}).




\renewcommand{\theequation}{\thesection.\arabic{equation}}
\setcounter{equation}{0}
\section{Preliminaries}\label{prelim}

In the following main result of this section we obtain bounds for the first and second moment of the permanent~\(T_n\) of the matrix~\(Y\) defined prior to Theorem~\ref{main_thm}. All constants throughout do not depend on~\(n.\)
\begin{Proposition} \label{t_n_sec_lem}
Let~\(r_{up},r_{low}\) and~\(\mu_n\) be as in Theorem~\ref{main_thm} and recall that~\(\delta = \mathbb{E}Z_{1,1}^2\) and~\(\nu = \mathbb{E}Z_{1,1}.\)
Defining
\begin{equation}\label{alp_def}
\alpha_{up} := \left(\frac{n(r_{up}-1)}{r_{up}(n-1)}\right)^{n} \text{ and } \beta_{up} := \frac{\delta r_{up}(n-1)}{\nu^2 r_{low}(r_{up}-1)}
\end{equation}
and
\begin{equation}\label{alp_def2}
\alpha_{low} := \left(\frac{n(r_{low}-1)}{r_{low}(n-1)}\right)^{n} \text{ and } \beta_{low} := \frac{\delta r_{low}(n-1)}{\nu^2 r_{up}(r_{low}-1)},
\end{equation}
we have that~\(\mathbb{E}T_n = \mu_n.\) Moreover, if~\(r_{low} \geq \frac{6\delta}{\nu^2},\) then we also have that
\begin{equation}\label{tn_bds}
\alpha_{low} \cdot e^{\beta_{low}-1} \cdot \left(1-\frac{2e}{n^2}\right)\leq \frac{\mathbb{E}T^2_{n}}{\mu_n^2} \leq \alpha_{up} \cdot e^{\beta_{up}-1} \cdot \left(1+\frac{2e}{n^2}\right)
\end{equation}
for all~\(n\) large.
\end{Proposition}
To illustrate, we again consider the example of~\(0-1\) matrices and set~\(Z_{i,j} = 1\) for all~\(i,j.\) If~\(r_i = r = r(n)\) for all~\(i\) so that the matrix is homogenous, then \[\alpha_{up} = \alpha_{low} = \left(\frac{n(r-1)}{r(n-1)}\right)^{n}  = e\left(1-\frac{1}{r}\right)^{n}(1+o(1))\] and~\(\beta_{up} = \beta_{low} = \frac{n-1}{r-1}.\) Plugging these into~(\ref{tn_bds}) and setting~\(\theta = \left(1-\frac{1}{r}\right) \cdot e^{\frac{1}{r-1}},\) we then get that \[\frac{\mathbb{E}T^2_{n}}{\mu_n^2} = e^{-\frac{1}{r-1}} \cdot \theta^{n} \cdot (1+o(1)).\] We remark that for constant~\(r,\) the term~\(\theta > 1\) strictly and therefore the second moment increases exponentially in~\(\theta.\) However, we show in the next Section that if~\(r = r(n)\) grows faster than~\(\sqrt{n},\) then the scaled variance of~\(T_n\) in fact converges to zero as~\(n \rightarrow \infty.\)



\subsection*{\em Proof of Proposition~\ref{t_n_sec_lem}}
The rows of~\(X\) are mutually independent and so for any~\(1 \leq i ,j \leq n\) we have that
\begin{equation}\label{p_dist}
\mathbb{P}(X_{i,j} = 1) = \frac{{n-1 \choose r_i-1}}{{n \choose r_i}} = \frac{r_i}{n} =:p_i
\end{equation}
and for any~\(1 \leq j_1\neq j_2 \leq n\) we also have that
\begin{equation}\label{p_dist2}
\mathbb{P}(X_{i,j_1} = 1,X_{i,j_2}=1) = \frac{{n-2 \choose r_i-2}}{{n \choose r_i}} = \frac{r_i(r_i-1)}{n(n-1)}.
\end{equation}

For a permutation~\(\sigma\) of~\(\{1,2,\ldots,n\},\) we define~\(R_{\sigma} := \prod_{i=1}^{n} X_{i,\sigma(i)}Y_{i,\sigma(i)}\)
and get from~(\ref{per_def}) that~\(T_n = per(Y) = \sum_{\sigma} R_{\sigma}\) where the summation  is over all permutations of~\(\{1,2,\ldots,n\}.\) For any~\(\sigma\) we get from~(\ref{p_dist}) that\\\(\mathbb{E}R_{\sigma} = \prod_{i=1}^{n} \nu p_i = \frac{\nu^{n}\prod_{i=1}^{n} r_i}{n^{n}}\) and so~\(\mathbb{E}T_n = \nu^{n} n!\frac{\prod_{i=1}^{n} r_i}{n^{n}} = \mu_n.\)

\emph{Proof of the upper bound in~(\ref{tn_bds})}: Expanding~\(T_n^2 = \left(\sum_{\sigma} R_{\sigma}\right)^2\) we get that
\begin{equation}\label{tn2}
\mathbb{E}T_n^2 = \sum_{\sigma_{1},\sigma_2} \mathbb{E}R_{\sigma_{1}} R_{\sigma_2} = \sum_{\sigma_1} \sum_{k=0}^{n} \sum_{{\cal S} : \#{\cal S} = k}\sum_{\sigma_2 : \sigma_1 \cap \sigma_2 = {\cal S}} \mathbb{E}R_{\sigma_1} R_{\sigma_2}
\end{equation}
where~\(\sigma_1 \cap \sigma_2 := \{1 \leq i \leq n : \sigma_1(i) = \sigma_2(i)\}\) is the set of indices where~\(\sigma_1\) and~\(\sigma_2\) agree. If~\(\sigma_1 \cap \sigma_2 = {\cal S}\) and~\(\#{\cal S} = k,\) then from~(\ref{p_dist2}) and the definition of~\(\mu_n,\) we see that
\begin{eqnarray}
\mathbb{E}R_{\sigma_1} R_{\sigma_2} &=& \prod_{j \in {\cal S}} \delta p_j \cdot \prod_{i \in {\cal S}^c}(\nu p_i)^2 \cdot \left(\left(1-\frac{1}{r_i}\right) \cdot \frac{n}{n-1}\right) \nonumber\\
&=& \left(\frac{\mu_n}{n!}\right)^2 \cdot \left(\frac{\delta}{\nu^2}\right)^{k} \cdot \frac{n^{n}}{(n-1)^{n-k}} \cdot \frac{1}{\prod_{j \in {\cal S}} r_j}  \cdot \prod_{i \in {\cal S}^{c}} \left(1-\frac{1}{r_i}\right). \nonumber\\
\label{r_first}
\end{eqnarray}
Using~\(r_{low} \leq r_i \leq r_{up}\) in the last two terms of~(\ref{r_first}) we then get that
\begin{eqnarray}
\mathbb{E}R_{\sigma_1} R_{\sigma_2} &\leq& \frac{\mu_n^2 \alpha}{(n!)^2} \cdot \beta^{k}\label{r_12_est_ax}
\end{eqnarray}
where~\(\alpha = \alpha_{up}\) and~\(\beta = \beta_{up}\) are as defined in~(\ref{alp_def}).

Substituting~(\ref{r_12_est_ax}) into~(\ref{tn2}) we get
\begin{equation}\label{del_4}
\mathbb{E}T_n^2 \leq \frac{\mu_n^2\alpha}{(n!)^2} \sum_{\sigma_1} \sum_{k=0}^{n} \sum_{{\cal S} : \#{\cal S} = k}\sum_{\sigma_2 : \sigma_1 \cap \sigma_2 = {\cal S}} \beta^{k}.
\end{equation}
Given~\(\sigma_1\) and~\({\cal S}\) with~\(\#{\cal S} = k,\) the number of permutations~\(\sigma_2\) satisfying\\\(\sigma_2 \cap \sigma_1 = {\cal S}\) is equal to the number of derangements of~\(\{1,2,\ldots,n-k\},\) which in turn equals~\((n-k)!\cdot b_{n-k},\) where~\(b_j := \sum_{l=0}^{j}\frac{(-1)^{l}}{l!}.\)
From~(\ref{del_4}) we therefore get that
\begin{eqnarray}
\mathbb{E}T_n^2  &\leq& \frac{\mu_n^2\alpha}{(n!)^2} \sum_{\sigma_1} \sum_{k=0}^{n} \sum_{{\cal S} : \#{\cal S} = k}(n-k)!\cdot b_{n-k}\cdot \beta^{k} \nonumber\\
&=& \frac{\mu_n^2\alpha}{(n!)^2} \sum_{\sigma_1} \sum_{k=0}^{n} {n \choose k} (n-k)!\cdot b_{n-k}\cdot \beta^{k} \nonumber\\
&=& \frac{\mu_n^2\alpha}{n!} \sum_{k=0}^{n} {n \choose k} (n-k)!\cdot b_{n-k}\cdot \beta^{k} \nonumber\\
&=& \mu_n^2\alpha \sum_{k=0}^{n}\frac{\beta^{k}}{k!}\cdot b_{n-k}.\label{del_5}
\end{eqnarray}

To estimate the final term in~(\ref{del_5}) we split the corresponding summation into two terms depending on whether~\(1 \leq k \leq n-\log{n}\) or~\(n-\log{n} \leq k \leq n.\) Using~\(|b_{n-k} -e^{-1}| \leq \sum_{j \geq n-k+1}\frac{1}{j!}\) we have for~\(0 \leq k \leq n-\log{n}\) that
\begin{equation}\label{bn_est}
|b_{n-k}-e^{-1}| \leq \sum_{ j\geq \log{n}} \frac{1}{j!} \leq \frac{1}{n^2}
\end{equation}
for all~\(n\) large and so
\begin{equation}\label{del_7}
\sum_{k=0}^{n-\log{n}}\frac{\beta^{k}}{k!} \cdot |b_{n-k}-e^{-1}| \leq \frac{1}{n^2}\sum_{k=0}^{n-\log{n}}\frac{\beta^{k}}{k!}.
\end{equation}
For the other range~\(n-\log{n} \leq k \leq n,\) we use~\(|b_{n-k}| \leq e\) to get that
\begin{equation}\label{del_6}
\sum_{k=n-\log{n}}^{n}\frac{\beta^{k}}{k!} \cdot |b_{n-k}-e^{-1}| \leq (e+e^{-1})\sum_{k=n-\log{n}}^{n} \frac{\beta^{k}}{k!}.
\end{equation}
Using~\(k! \geq k^{k} e^{-k}\) and~\(\beta \leq \frac{2n\delta}{r_{low}\nu^2}\)  we get that the final term in~(\ref{del_6}) is at most
\[ (e+e^{-1})\sum_{k=n-\log{n}}^{n}\left(\frac{2ne}{kr_{low}}\right)^{k} \leq (e+e^{-1})\sum_{k=n-\log{n}}^{n} \left(\frac{2ne \delta}{(n-\log{n})r_{low}\nu^2}\right)^{k}.\]
Since~\(r_{low} \geq \frac{6\delta}{\nu^2}\) (see statement of Proposition~\ref{t_n_sec_lem}) we have~\(\frac{2ne\delta}{(n-\log{n})r_{low}\nu^2} \leq e^{-2C}\) for all~\(n\) large and some positive constant~\(C.\) Thus
\begin{equation}\label{del_9}
\sum_{k=n-\log{n}}^{n} \left(\frac{2ne \delta}{(n-\log{n})r_{low}\nu^2}\right)^{k} \leq \sum_{k \geq n-\log{n}}e^{-2Ck} \leq e^{- Cn}
\end{equation}
for all~\(n\) large.

From~(\ref{del_7}),~(\ref{del_6}),~(\ref{del_9}) and the above discussion we then get that
\begin{eqnarray}
\sum_{k=0}^{n}\frac{\beta^{k}}{k!} \cdot b_{n-k} &\leq& \left(e^{-1}+\frac{1}{n^2}\right)\sum_{k=0}^{n}\frac{\beta^{k}}{k!} +(e+e^{-1})e^{-Cn} \nonumber\\
&\leq& \left(e^{-1}+\frac{1}{n^2}\right)e^{\beta} + (e+e^{-1})e^{-Cn} \nonumber\\
&\leq& \left(e^{-1}+\frac{2}{n^2}\right)e^{\beta} \label{del_8}
\end{eqnarray}
for all~\(n\) large. Substituting~(\ref{del_8}) into~(\ref{del_5}) and using~(\ref{tn2}), we obtain the upper bound in~(\ref{tn_bds}).~\(\qed\)

\emph{Proof of the lower bound in~(\ref{tn_bds})}: Using~(\ref{r_first}) and arguing as in the paragraph preceding~(\ref{del_5})  we have that
\begin{equation}\label{del_low2}
\mathbb{E}T_n^2  \geq \mu_n^2 \cdot \alpha_{low} \sum_{k=0}^{n}\frac{\beta^{k}_{low}}{k!}\cdot b_{n-k}
\end{equation}
and using~(\ref{bn_est}) we get that the term~\(\sum_{k=0}^{n}\frac{\beta^k_{low}}{k!} \cdot b_{n-k} \) is bounded below by
\begin{equation}\label{del_low}
\left(e^{-1}-\frac{1}{n^2}\right)\sum_{k=0}^{n-\log{n}}\frac{\beta^{k}_{low}}{k!}
\end{equation}
for all~\(n\) large. Using~\(r_{up} \geq r_{low} \geq \frac{6\delta}{\nu^2}\) and arguing as in~(\ref{del_6}) and~(\ref{del_9}), we get that
\begin{equation}\label{del_low3}
\sum_{k \geq n-\log{n}} \frac{\beta^{k}_{low}}{k!}\leq e^{-D n}
\end{equation}
for some constant~\(D > 0\) and all~\(n\) large. From~(\ref{del_low2}),~(\ref{del_low}) and~(\ref{del_low3}), we then get that
\begin{eqnarray}
\mathbb{E}T_n^2  &\geq& \mu_n^2 \alpha_{low} \left(\left(e^{-1}-\frac{1}{n^2}\right)e^{\beta_{low}} - \left(e^{-1}-\frac{1}{n^2}\right)e^{-Dn} \right)\nonumber\\
&\geq& \mu_n^2 \cdot \alpha_{low}\cdot e^{\beta_{low}-1}\left(1-\frac{2e}{n^2}\right) \nonumber
\end{eqnarray}
for all~\(n\) large. Substituting this into~(\ref{tn2}), we get the lower bound in~(\ref{tn_bds}).~\(\qed\)


\renewcommand{\theequation}{\thesection.\arabic{equation}}
\setcounter{equation}{0}
\section{Proof of Theorem~\ref{main_thm}}\label{main_sec}
As before~\(T_n\) is the permanent of the constrained random matrix~\(Y.\) We first find an upper bound for the variance of~\(T_n\) using Proposition~\ref{t_n_sec_lem}. Let~\(r_{up} = r\) and let~\(\alpha_{up}\) and~\(\beta_{up}\) be as in~(\ref{alp_def}). Using~\(\left(1-\frac{1}{n}\right)^{-n}  = e(1+o(1))\) where~\(o(1) \longrightarrow 0\) as~\(n\rightarrow \infty,\) we then have that~\[e^{-1} \cdot \alpha_{up} = e^{-1} \cdot \left(1-\frac{1}{n}\right)^{-n} \cdot \left(1-\frac{1}{r}\right)^{n} = \left(1-\frac{1}{r}\right)^{n}(1+o(1))\] and~\(\beta_{up} = \frac{\delta r(n-1)}{\nu^2 r_{low}(r-1)} \leq \frac{\delta rn}{\nu^2 r_{low}(r-1)}.\) Thus
\begin{equation}\label{alp_eval}
\alpha_{up} \cdot e^{\beta_{up}-1}  = e^{-1} \alpha_{up} \cdot e^{\beta_{up}} \leq \left(e^{\frac{\delta r}{\nu^2 r_{low}(r-1)}} \cdot \left(1-\frac{1}{r}\right)\right)^{n} \cdot (1+o(1)).
\end{equation}




For~\(r \geq 2\) we have that~\(\frac{r}{r-1} \leq 2\) and moreover, from~(\ref{r_cond}) we know that~\(\frac{\delta}{\nu^2 r_{low}} \leq \frac{1}{4}\) for all~\(n\) large. This implies that~\(\frac{\delta r}{\nu^2 r_{low}(r-1)} \leq \frac{2\delta}{\nu^2r_{low}} \leq \frac{1}{2}\) for all~\(n\) large and so using standard~\(O\) notations we have that
\begin{eqnarray}
e^{\frac{\delta r}{\nu^2 r_{low}(r-1)}} &=& 1+\frac{\delta r}{\nu^2 r_{low}(r-1)} + O\left(\frac{\delta}{\nu^2 r_{low}}\right)^2 \nonumber\\
&=& 1+\frac{\delta}{\nu^2r_{low}} + \frac{\delta}{\nu^2r_{low}(r-1)} + O\left(\frac{\delta}{\nu^2 r_{low}}\right)^2 \nonumber\\
&=& 1 + \frac{\delta}{\nu^2 r_{low}} + O\left(\frac{\delta}{\nu^2 r_{low}}\right)^2 \label{plug_b}
\end{eqnarray}
since~\(r = r_{up} \geq r_{low}\) and~\(\delta = \mathbb{E}Z^2 \geq (\mathbb{E}Z)^2 = \nu^2.\) Thus
\begin{eqnarray}
e^{\frac{\delta r}{\nu^2 r_{low}(r-1)}} \cdot \left(1-\frac{1}{r}\right)
&=& \left(1+\frac{\delta }{\nu^2 r_{low}} + O\left(\frac{\delta}{\nu^2 r_{low}}\right)^2\right) \cdot \left(1-\frac{1}{r_{up}}\right) \nonumber\\
&=& 1+\frac{\delta}{\nu^2 r_{low}} - \frac{1}{r_{up}} -\frac{\delta }{\nu^2 r_{low}r_{up}} +O\left(\frac{\delta}{\nu^2 r_{low}}\right)^2\nonumber\\
&=&1 +\frac{\delta}{\nu^2 r_{low}} - \frac{1}{r_{up}} + O\left(\frac{\delta}{\nu^2 r_{low}}\right)^2\nonumber\\
&=& 1+ o\left(\frac{1}{n}\right)\label{plug_a}
\end{eqnarray}
where as before, the third expression in~(\ref{plug_a}) follows from the fact that\\\(r_{up} \geq r_{low}\) and~\(\delta \geq \nu^2\) and the final expression in~(\ref{plug_a}) follows from the conditions in~(\ref{r_cond}).

Plugging~(\ref{plug_a}) into~(\ref{alp_eval}) we get that
\begin{equation}
\alpha_{up} \cdot e^{\beta_{up}-1} \leq \left(1+o\left(\frac{1}{n}\right)\right)^{n} \cdot (1+o(1))= 1+o(1). \nonumber
\end{equation}
From~(\ref{tn_bds}) we therefore get that~\(\frac{\mathbb{E}T_n^2}{\mu_n^2} \leq 1+o(1)\) or equivalently that\\\(var(T_n) \leq o(\mu_n^2).\) This obtains~(\ref{m_conv}) and therefore completes the proof of Theorem~\ref{main_thm}.~\(\qed\)

\section*{Acknowledgement}
I thank Professors Rahul Roy, C. R. Subramanian and the referee for crucial comments that led to an improvement of the paper. I also thank IMSc for my fellowships.


\bibliographystyle{plain}

\end{document}